\documentclass[10pt]{article}
\expandafter\let\csname equation*\endcsname\relax
\expandafter\let\csname endequation*\endcsname\relax
\usepackage{bm}
\usepackage{enumerate}
\usepackage[top=1in, bottom=1in, left=1.2in, right=1.2in]{geometry}
\usepackage{mathrsfs,color}
\usepackage{amsfonts}\usepackage{multirow}
\usepackage{amssymb,dsfont}
\usepackage{commath}
\usepackage{caption,comment}
\usepackage{blkarray}
\usepackage{amsmath, amscd}
\usepackage{tikz-cd}
\usepackage{float}
\usepackage{amssymb,amsthm}
\usepackage[toc,page]{appendix}
\usepackage{graphicx,afterpage}
\usepackage{epstopdf}
\usepackage{chngcntr}
\usepackage{apptools}
\AtAppendix{\counterwithin{theorem}{section}}

\usepackage[pdftex,colorlinks=true]{hyperref}
\hypersetup{CJKbookmarks,%
	bookmarksnumbered,%
	colorlinks,%
	linkcolor=black,%
	citecolor=black,%
	plainpages=false,%
	pdfstartview=FitH}
\numberwithin{equation}{section}
\numberwithin{figure}{section}

\newcommand\tabcaption{\def\@captype{table}\caption}

\newtheorem{thm}{Theorem}[section]

\newtheorem{lem}[thm]{Lemma}

\newtheorem{aspt}[thm]{Assumption}
\newtheorem{rem}[thm]{Remark}

\newtheorem{theorem}{Theorem}[section]

\allowdisplaybreaks

\newcommand{\normo}[1]{\left\lvert#1\right\rvert}

 \usepackage{array}
 \newcolumntype{M}[1]{>{\centering\arraybackslash}m{#1}}
 \newcolumntype{N}{@{}m{0pt}@{}}
\usepackage[toc,page]{appendix}

\usepackage{authblk}

\date{}

\begin{document}
	
	\title{$L^2$ convergence of smooth approximations of Stochastic Differential Equations with unbounded coefficients}
	
	\author[ ]{Sahani Pathiraja}
	\affil[ ]{ \small{Institute of Mathematics, University of Potsdam, Germany}}
	
	\maketitle

		\begin{abstract} 
		
		The aim of this paper is to obtain convergence in mean in the uniform topology of piecewise linear approximations of Stochastic Differential Equations (SDEs) with $C^1$ drift and $C^2$ diffusion coefficients with uniformly bounded derivatives.  Convergence analyses for such Wong-Zakai approximations most often assume that the coefficients of the SDE are uniformly bounded.  Almost sure convergence in the unbounded case can be obtained using now standard rough path techniques, although $L^q$ convergence appears yet to be established and is of importance for several applications involving Monte-Carlo approximations.  We consider $L^2$ convergence in the unbounded case using a combination of traditional stochastic analysis and rough path techniques.  We expect our proof technique extend to more general piecewise smooth approximations. 
		
	\end{abstract}
	\noindent
	{\bf Keywords.} Wong-Zakai, unbounded coefficients, piecewise smooth approximations, stochastic differential equations, rough paths \\

	\vspace{0.5cm}

	\section{Introduction}
	\label{sec:intro}
	
	Given a filtered probability space $(\Omega, \mathcal{F}, \mathbb{P}, \{\mathcal{F}_{t} \}_{t \in [0,T]})$ satisfying the usual conditions and a time interval $[0,T], T > 0$, we are interested in a stochastic process $\{X_t\}_{t \in [0,T]}$ given by 
	\begin{align}
	\label{eq:origsde}
	X_t = x + \int_0^t b(X_s)ds + \int_0^t \sigma(X_s) dB_s, 
	\end{align}
	where $X_t \in \mathbb{R}^v$, $b: \mathbb{R}^v \rightarrow \mathbb{R}^v$ and $\sigma: \mathbb{R}^v \rightarrow \mathbb{R}^{v \times r}$ satisfy appropriate regularity conditions and $B$ is a $r$-dimensional Wiener process.  The stochastic integrals in the rhs of (\ref{eq:origsde}) are interpreted in the Ito sense.  It is often necessary in many applications to utilise smooth approximations of the Brownian motion and to consider the corresponding solutions to (\ref{eq:origsde}). To this end, let $\Pi_d$ denote a partition of the time interval $[0,T]$ consisting of $d$ units, $d \in \mathbb{N}$ and associated filtration $\mathcal{F}^d := \sigma(B_t : t \in \Pi_d)$.  In the case where a sequence of approximations $(B^d)_{d \geq 1}$ are $\mathcal{F}^d$-martingales and satisfy a uniform tightness condition \cite{Jakubowski1989}, it is known that (\ref{eq:origsde}) with $dB_s^d$ instead of $dB_s$ converges in distribution to (\ref{eq:origsde})\footnote{This holds for more general semimartingale drivers, but for the purposes of this paper, we focus on the case of Wiener processes.} \cite{Kurtz1991}.  Many piecewise smooth approximations do not satisfy these conditions so that the limiting process involves additional correction terms.  More precisely, consider $B^d(\omega)$ a piecewise smooth approximation of the Brownian motion $B(\omega)$, with $B$ as its limiting process, and the ODE 
	\begin{align}
	\label{eq:origsdedelt}
	dX^d_t = b(X^d_t)dt + \sigma(X^d_t) \dot{B}_t^d dt.
	\end{align}
	In the above and throughout the article, we drop the $\omega$ for notational ease.  It is now well-known that under certain regularity conditions on the coefficients and conditions on the construction of $B_t^d$,  solutions of (\ref{eq:origsdedelt}) do not converge to solutions of (\ref{eq:origsde}) as the mesh size $\delta \rightarrow 0$, but rather to a Stratonovich SDE plus an anti-symmetric drift term that depends on how the approximation is constructed, i.e. 
	\begin{align}
	\label{eq:stratantisde}
	dX_t = b(X_t)dt + \sigma(X_t) \circ dB_t + \sum_{n,m = 1}^v S_{n,m}[A_nX_t,A_mX_t] dt, 
	\end{align}
	where the notation $C_{i,j}$ denotes the $i,j$-th entry of a matrix $C$, $[\cdot, \cdot]$ denotes the Lie bracket and the operator $A_n$ is defined by $A_nx := \sum_i \sigma_{i,n}(x) \frac{\partial}{\partial x_i} \enskip \forall \enskip n = 1, 2, \cdots r$.  One important condition for the previous statement is the existence of a skew symmetric matrix $S$ as the limit of the area process 
	\begin{align*}
	\frac{1}{2 \delta}\mathbb{E} \left[ \int_0^\delta B^d_s \otimes \dot{B}_s^d -  \dot{B}_s^d \otimes B_s^d ds \right].
	\end{align*}
	For reference, see Theorem 7.2 in \cite{Ikeda1989} which captures much of the pioneering work on this matter for the multi-dimensional case by e.g. \cite{Sussmann1991}, \cite{McShane1972}, \cite{Stroock1972}, \cite{Nakao1976}, \cite{Ikeda1977}.  Such piecewise smooth approximation models are relevant in many areas of mathematical modelling where it is necessary to work with ODEs with random fluctuations, e.g. for extracting macroscopic models from detailed microscopic stochastic dynamics (e.g. \cite{Huisinga2003}, \cite{Pavliotis2008}).  Secondly such approximations are useful because they are naturally tied to the Stratonovich interpretation of stochastic integrals, in particular, when the anti-symmetric term in (\ref{eq:stratantisde}) is zero, which holds in the scalar case, when the vector fields $A_n$ commute, or when $S \equiv 0$ which occurs for piecewise linear approximations or mollifiers \cite{Ikeda1989}.  Continuity of the solutions of ODEs driven by such smooth approximations with respect to the driving path in the uniform topology is ensured by the Doss-Sussman theorem.  This property is advantageous for various applications, in particular for robust filtering \cite{Clark2005}.  From a numerical approximation point of view, it is possible to utilise the many techniques available for ODEs (rather than SDEs) to then construct time-discrete approximations to (\ref{eq:stratantisde}).  This is in constrast to more direct numerical approximation schemes applied to Ito SDEs (\ref{eq:origsde}) like the Euler-Maruyama or Milstein which converge strongly to (\ref{eq:origsde}) as the time step tends to zero, i.e. without any additional correction term.  The convergence analysis of Wong-Zakai approximations is complicated by the presence of anticipative integrands, so that the traditional martingale techniques cannot be directly applied as in the analysis of Euler-Maruyama and related schemes.

The literature on convergence analyses of piecewise smooth approximations for both SDEs and SPDEs (e.g. \cite{HairerPardoux2015}, \cite{Gyongy2006}, \cite{Diehl1991}, \cite{Tessitore}, \cite{Hu2002}, \cite{Brzezniak1995}) is enormous, starting from the seminal work of \cite{Eugene1965} where convergence in probability of (\ref{eq:origsdedelt}) to (\ref{eq:stratantisde}) with $a = 0$ (antisymmetric term) in the scalar case was established.  Much of the inital work in this direction for SDEs relies on standard stochastic analysis techniques, loosely speaking, by introducing time shifts so that martingale inequalities can be applied (e.g. \cite{Gyongy2004}, \cite{Ikeda1989}).  In these works, uniform boundedness of the coefficients is crucial in order to control the extra terms introduced by the time shift.  More recently, \cite{Liu2019} showed convergence in probability in the Skorokhod $M1$ topology of Wong-Zakai approximations of scalar SDEs driven by a general semimartingale (thereby extending the earlier work of \cite{Protter1985} and \cite{Kurtz1995}), again for uniformly bounded coefficients.   Over the last couple of decades, the theory of rough paths pioneered by T. Lyons \cite{Lyons1998} has been utilised to obtain convergence results for a broader class of drivers and smooth approximations (see e.g. \cite{Frizober2009}, \cite{Lejay2006}, \cite{Coutin2007}, \cite{Friz2014}).  After defining an appropriate lift of the driving path (e.g. Stratonovich lift in the case of piecewise linear approximations), and by ensuring that the pathwise solution of the rough differential equation (RDE) coincides with the limiting SDE, convergence is easily established using continuity of the Ito-Lyons map (i.e. the solution map corresponding to the RDE, as a function of the initial condition and driver) (see e.g. Theorem 9.3 in \cite{Friz2014}).  Generally, weak convergence of smooth approximations of the driver then imply weak convergence of the (random) ODE solutions.  Strong convergence results can be obtained in specific situations; e.g. when $b$ is $Lip(1+\epsilon)$ and $\sigma$ is $Lip(2 + \epsilon)$ in the sense of Stein, for dyadic piecewise linear approximations of Brownian motion \cite{Coutin2007} (see also \cite{Friz2014} for the non-dyadic case).  A common feature of these works is the lipschitz assumption as in the Stein sense \cite{Frizober2009}, \cite{Coutin2007}, therefore requiring uniform boundedness of the coefficients.  \cite{KellyMelb2016} use rough path techniques to obtain weak convergence under much weaker conditions ($C^{1+\epsilon}$ and $C^{2+\epsilon}$ for the drift and diffusion coefficients respectively) and for a broad range of smooth approximations satisfying a weak invariance principle.  As noted in \cite{Friz2014}, it is in general not straightforward to use rough path techniques for $L^q$ convergence, although Chapter 10 of \cite{Friz2014} and also \cite{Friza2010} presents some analysis in this direction in the case of Gaussian processes.  To summarise, the aforementioned results apply to the case of bounded coefficients and/or weak convergence or almost sure convergence.
	
Our aim is to obtain $L^2$ rates of convergence of piecewise smooth approximations of (\ref{eq:origsde}) with unbounded coefficients.  $L^q$ convergence of numerical approximations of SDEs is often of central importance in many Monte Carlo sampling methods such as in Markov Chain Monte Carlo \cite{Hutzenthaler2012} and in control-type sequential Monte Carlo methods which rely on empirical approximations of McKean-Vlasov SDEs \cite{P_mckean_arxiv2020}.   
Here we present a simple approach that combines elements of stochastic analysis and rough path techniques to achieve this result, specifically, focusing on nested piecewise linear approximations of the Brownian path.  That is, $\Pi_d$ is a partition of the time interval $[0,T]$ with associated mesh size $\delta_d$ such that $\Pi_d \subset \Pi_{d+1} \enskip \forall \enskip d \in \mathbb{N}$ and $\delta_d \rightarrow 0$ as $d \rightarrow \infty$ and 
	\begin{align*}
	B_t^d  = B_{t_i}   + \frac{t - t_i}{\delta_d} (B_{t_{i+1}} - B_{t_i}), \quad \text{for} \enskip t \in [t_i, t_{i+1}). 
	\end{align*}
	 We expect that our result can be extended to more general approximations using the same proof structure.  
	
	\subsection{Statement of Main Result}
	Before stating our main result, we state the main assumptions utilised throughout the paper.

	\begin{aspt}
			\label{ass:bsig}
		 $b \in C^1(\mathbb{R}^v, \mathbb{R}^v)$ and $\sigma \in C^2(\mathbb{R}^v, \mathbb{R}^{v \times r})$  and there exists some constant $M_c$ such that 
		\begin{align*}
		\sup_{x \in \mathbb{R}^v, \enskip j \in \{1, 2, \cdots r \}} \normo{\nabla b} \vee \normo{\nabla \sigma_{\cdot, j}}  \vee \normo{\nabla^2 \sigma_{\cdot, j}} \leq M_c 
		\end{align*}
		where $\sigma_{\cdot, j}$ denotes the $j$th column of the matrix $\sigma$.
			\end{aspt}

\begin{aspt}
		\label{ass:correc}
		The vector function $\Sigma: \mathbb{R}^v \rightarrow \mathbb{R}^v$ with $i$-th entry given by   
		\begin{align*}
		\Sigma_i (x) := \sum_{k,m=1}^r \sum_{l=1}^v \sigma_{l,k}(x) \frac{\partial }{\partial x^l} \sigma_{i,m}(x)
		\end{align*}
		is globally Lipschitz-continuous.
	\end{aspt}

	These assumptions ensure well-posedness of the ODE and limiting SDE, as well as the RDE corresponding to the Stratonovich lift of $B^d$.  Assumption \ref{ass:correc} in particular plays a crucial role in obtaining a pathwise control on solutions to the RDE, see Lemma \ref{lem:uniformest_sde}.

	\begin{theorem}
		\label{theo: sdeconv}
		Assume (\ref{ass:bsig}) and (\ref{ass:correc}).  Suppose $X_t^d$ is the unique solution of 
		\begin{align}
		\label{eq:smoothsde}
		dX_t^d = b(X_t^d)dt + \sigma(X_t^d) \dot{B}_t^d dt 
		\end{align}
		and $X_t$ is the unique strong solution of
		\begin{align}
		\label{eq:stratsde}
		dX_t = b(X_t)dt + \sigma(X_t) \circ dB_t. 
		\end{align}
		with $X_0^d = X_0 = x_0$ for some fixed but arbitrary $x_0$.  It holds that for every $T > 0$, 
		\begin{align*}
		\lim_{d \rightarrow \infty} \mathbb{E} \left[ \sup_{0\leq t \leq T} \normo{X_t - X_t^d}^2  \right] = 0.
		\end{align*}
		
	\end{theorem}

	\begin{rem}
		Almost sure convergence under Assumptions \ref{ass:bsig} and \ref{ass:correc} can be obtained using continuity of the Ito-Lyons map established in \cite{Lejayglob2012} and standard rough path arguments. 
	\end{rem}
	
	\section{Notation and Background on rough path theory}
	\label{sec:notandroughpath}
	
	We briefly summarise some of the necessary background on rough path theory that will be utilised in Lemma \ref{lem:uniformest_sde}.  This summary is by no means exhaustive and much of the algebraic and geometric technicalties are left out as they are not needed for understanding the main results. We refer to the extensive literature on rough paths (e.g. \cite{Friz2014}, \cite{Friz2009}, \cite{Lyons2002}) for further details.  
	
	A rough path framework allows for developing pathwise solutions to SDEs by considering a fixed realisation of the driver (i.e. $B(\omega)$ for fixed $\omega$, in our case), so that the differential equation is seen rather as a controlled ODE.  The irregularity of the driver means that special care is needed to define solutions to such ODEs.  The fundamental insight from rough path theory is that   $n$-order iterated integrals of the form\footnote{the notation $p_i$ is used to refer to the $i$th index of a permutation of the set $\{1, 2, \cdots r \}$}  $\int_s^t dX_{t_1}^{p_1} \int_s^{t_1} dX_{t_2}^{p_2} \cdots \int_s^{t_{n-1}} dX_{t_n}^{p_n}$ with $n \leq \lfloor{\frac{1}{\alpha}} \rfloor$ for a driver $X$ of $\alpha$-Hoelder regularity (and an apriori interpretation of these iterated integrals) is sufficient to make these differential equations well-defined in a deterministic sense.  The iterated integrals or ``signature'' can be seen heuristically as arising from a Taylor expansion \cite{Lejayintro2003}.  Since Brownian paths are $\alpha$-Hoelder regular with $\alpha \in (1/3, 1/2)$ it is only necessary to consider the second order iterated integral taking the form $\int_s^t (B_{t_1} - B_s) \otimes dB_{t_1} \in \mathbb{R}^r \otimes \mathbb{R}^r$.  There is still some ambiguity as to how this stochastic integral should be interpreted; a Stratonovich (Ito) interpretation corresponds to the so-called Stratonovich (Ito)-lift or enhancement.  We use the notation ${\bf X} := ({\bf X}^1, {\bf X}^2)$ to denote the rough path lift of a continuous path $X: [0,T] \rightarrow \mathbb{R}^r$, where ${\bf X}^1 := X$ and ${\bf X}^2: [0,T]^2 \rightarrow \mathbb{R}^r \otimes \mathbb{R}^r$ denotes the second order process such that ${\bf X}^2_{s,t} =:\int_s^t X_{s,t_1} \otimes dX_{t_1}$ where $X_{s,t} := X_t - X_s$.  We are concerned with Brownian rough paths; these paths belong to the space of $\alpha$-Hoelder rough paths $C^\alpha ([0,T], \mathbb{R}^r)$ such that 
	\begin{align*}
	\norm{{\bf X}}_\alpha := \sup_{0 \leq s < t \leq T} \max \left\{  \frac{\normo{{\bf X}_{s,t}^1}}{|t-s|^\alpha}  ,\frac{\normo{{\bf X}_{s,t}^2}}{|t-s|^{2 \alpha}} \right\} 
	\end{align*}
	is finite and Chen's relation is satisfied.  For an $\alpha$-Hoelder continuous path $X$ taking values in $\mathbb{R}^m$, we define 
	\begin{align*}
	\norm{X}_\alpha := \sup_{0 \leq s < t \leq T} \frac{\normo{X_{s,t}}}{|t-s|^\alpha}.
	\end{align*}
	Throughout the article, we denote by ${\bf B}$ the Stratonovich lift of $B$, which is a geometric rough path.  These have the important property that there exists a sequence of smooth (or piecewise smooth) rough paths converging uniformly on $[0,T]$ in the rough path topology to the lifted path (see e.g. Proposition 2.5 in \cite{Friz2014}).  Likewise, we denote by ${\bf B}^d$ the Stratonovich lift of a piecewise smooth approximation $B^d$.      
	  
	It is classical that the solution of (\ref{eq:origsdedelt}) coincides with the solution (when it exists and is unique) of the following Rough Differential Equation (RDE) driven by ${\bf B}^d$  
	\begin{align}
	\label{eq:smoothrde}
	\hat{X}_t^d = \hat{X}_0 + \int_0^t b(\hat{X}_t^d)dt + \int_0^t \sigma(\hat{X}_t^d) {\bf B}_t^d,
	\end{align}
	(see e.g. Section 9 in \cite{Lejayglob2012}).  Roughly speaking, the integrand must be ``smooth enough'' to counteract the irregularity of the driver \cite{Unterberger2012}.  Several concepts of solutions to RDEs have been developed, e.g. solution in the sense of Lyons \cite{Lyons1998}, Gubinelli's controlled paths formulation \cite{Gubinelli2004}, Davie's formulation in terms of a Taylor expansion or Euler discretisation \cite{Davie2008} and also more recently the coordinate independent definition from \cite{Bailleul2015}.  There exist equivalencies between various solution concepts \cite{Bailleul2018},  e.g. solutions in the sense of Lyons and Davie are equivalent when $\sigma$ is $\gamma$-Hoelder continuous with $\alpha (2 + \gamma) > 1$, see \cite{Lejayglob2012}.  Here we consider solutions in the sense of Davie \cite{Davie2010} i.e. where $\hat{X}^d$ is a continuous path from $[0,T] \rightarrow \mathbb{R}^v$ of finite $p$-variation such that for some constant $L$,
	\begin{align}
	\label{eq:solndavie}
	\left|\hat{X}^d_t - \hat{X}^d_s - \int_s^t b(\hat{X}_u) du - \sigma(\hat{X}_s){\bf B}_{s,t}^{d,1} - \Sigma(\hat{X}_s){\bf B}_{s,t}^{d,2} \right| \leq L |t-s|^\theta		\quad \forall \enskip 0 \leq s < t \leq T,
	\end{align}
	where $\theta:= 3\alpha$.  With a slight abuse of notation, we denote by $\sigma$ a linear vector field from $\mathbb{R}^v \rightarrow \mathcal{L}(\mathbb{R}^r, \mathbb{R}^v)$ and $\Sigma(x)$ is also a linear vector field from $\mathbb{R}^v \rightarrow \mathcal{L}(\mathbb{R}^r \otimes \mathbb{R}^r, \mathbb{R}^v)$.

	\section{Proof of Theorem \ref{theo: sdeconv}}
	The proof involves utilising localisation arguments to extend well-known convergence results in the case of $C_b^2$ coefficients (i.e. the space of twice continuously differentiable uniformly bounded functions with bounded derivatives) as in \cite{Ikeda1989} to unbounded coefficients.  The main goal is to close the commutative diagram in Figure \ref{fig:figure1}.  We do so by establishing uniform in the approximation parameter $d$ and pointwise in the localisation parameter $n$ convergence results which thereby permit the interchange of limits in $n$ and $d$ using the Moore-Osgood theorem.  We also develop uniform in $d$ pathwise moment bounds on $X^d$ in Lemma \ref{lem:uniform} using rough path techniques, which plays a crucial role in establishing uniform convergence. 
	
	\begin{figure}[h]
	\begin{center}

$\begin{CD}
X^{d,n} @> d \rightarrow \infty > > X^n \\
@V n \rightarrow \infty V V @V V n \rightarrow \infty V\\
X^d @> d \rightarrow \infty > ? > X
\end{CD}$
	\end{center}
	\caption{Limit diagram. Convergence is always $L^2$ with respect to the uniform metric over $[0,T]$}
	\label{fig:figure1}
\end{figure}

	We start by considering a localised form of the limiting SDE in Ito form.  For all $n \in \mathbb{N}$, let $S_n := \{x \in \mathbb{R}^v \enskip | \enskip  |x|^2 \leq n \}$.  Consider $X_t^n$ which is the unique strong solution of 
	\begin{align}
	\label{eq:sden}
	dX_t^n = b^n(X_t^n)dt + \Sigma^n (X_t^n)dt + \sigma^n(X_t^n)dB_t 
	\end{align}
	where 
	\begin{align*}
	b^n(x) =
	\begin{cases}
	b(x), \enskip x \in S_n  \\
	f_b(x), \enskip x \notin S_n 
	\end{cases} 
	\end{align*}
	where $f_b \in C_b^1(\mathbb{R}^v, \mathbb{R}^v)$ 
	is chosen such that $b^n \in C_b^1(\mathbb{R}^v, \mathbb{R}^v)$.  The functions $\sigma^n$ and $\Sigma^n$ are constructed in a similar manner.
	Consider also the following ODE utilising the smooth approximation of the Brownian motion as described above
	\begin{align}
	\label{eq:smoothlocapp}
	dX_t^{d, n} = b^n(X_t^{d, n})dt + \sigma^n(X_t^{d, n}) \dot{B}_t^d dt.
	\end{align}
	We will assume throughout that $X_0^{n} = X_0^{d,n} = x_0$ where $x_0$ is fixed but arbitrary.  The following lemma utilises standard results to establish convergence along the clockwise route starting from $X^{d,n}$ in Figure \ref{fig:figure1}. 
	\begin{lem}
		\label{lem:locsde}
		For every $T > 0$ and $X_t^{d,n}$ satisfying (\ref{eq:sden}) and $X_t$ satisfying (\ref{eq:stratsde}), 
		\begin{align}
		\label{eq:firstdoub}
		\lim_{n \rightarrow \infty} \lim_{d \rightarrow \infty} \mathbb{E} \left[\sup_{0 \leq t \leq T} |X_t^{d,n} - X_t^n|^2  + \sup_{0 \leq t \leq T} |X_t^{n} - X_t|^2 \right]= 0.
		\end{align}
		
		\begin{proof}
		Since the coefficients of (\ref{eq:sden}) and (\ref{eq:stratsde}) are continuously differentiably and uniformly bounded with bounded derivatives, we have from Theorem 7.2 in \cite{Ikeda1989} that for every $T > 0$,	
		\begin{align}
		\label{eq:deltan}
		\lim_{d \rightarrow \infty} \mathbb{E} \left[\sup_{0 \leq t \leq T} |X_t^{d,n} - X_t^n|^2 \right] = 0, \quad \forall \enskip n \in \mathbb{N}.  
		\end{align}  		
		Since the coefficients $b^n, \sigma^n, b, \sigma$ are globally Lipschitz continuous, we have from standard localisation arguments that 
		\begin{align*}
		\lim_{n \rightarrow \infty} \mathbb{E} \left[\sup_{0 \leq t \leq T} |X_t^{n} - X_t|^2 \right] 
		\end{align*}	
		which gives the desired result. 

		\end{proof}
	\end{lem}

	As discussed in Section \ref{sec:notandroughpath}, the solution of the RDE (\ref{eq:smoothrde}) coincides with $X_t^d$, so that we can now work with the rich toolbox offered by rough path theory to obtain a pathwise moment bound on $X^d$, uniformly in $d$ (Lemma \ref{lem:uniformest_sde}).  An analogous result to Lemma \ref{lem:uniformest_sde} using stochastic analysis techniques in the case of bounded coefficients can be found in Lemma 7.2 in \cite{Ikeda1989}, where the now standard approach of applying a time shift to obtain martingales has proven difficult to adapt to the case of unbounded coefficients.  The rough path framework allows to side-step this difficulty, particularly since dealing with RDE solutions in the sense of Davie means we can avoid having to directly control the integrals driven by $B_t^d$.   
	
	Firstly, existence and uniqueness of solutions (in the sense of Davie, and equivalently in the sense of Lyons) to (\ref{eq:smoothrde}) under Assumption \ref{ass:bsig}-\ref{ass:correc} follows in a straightforward manner from Proposition 3 and Theorem 1 for driftless RDEs in \cite{Lejayglob2012}.  We obtain the following bound on the solution to (\ref{eq:smoothrde}) in the small time horizon $\mu:= T^\alpha \leq K^*$ by a straightforward extension of the boundedness of solutions to driftless RDEs in Proposition 2 of \cite{Lejayglob2012}, where the crucial point is to utilise the sewing lemma to obtain an expression for $L$.  The proof of the following Lemma can be found in the Appendix.  

	\begin{lem}
			 	\label{lem:ybound} \textbf{Boundedness of solution to (\ref{eq:smoothrde}) in the short time horizon.}
		Under Assumptions \ref{ass:bsig} and \ref{ass:correc}, and with $\alpha \in (1/3,1/2)$ for Brownian motion and recalling the notation $\mu:= T^\alpha$, $\theta:= 3\alpha$, we have that 
		\begin{align}
		\label{eq:yroughbound}
		\norm{\hat{X}^d}_\alpha  \leq C_{15}(\mu)  + C_{16}(\mu) |x_0|,
		\end{align}
		where 
		\begin{align*}
		C_{15} (\mu) &:= C_{13} \left(  C_{14}\mu^{2} \left( \norm{\nabla \Sigma}_\infty \left((1 + \mu^{2})\norm{{\bf B}^d}_\alpha^2 + \norm{{\bf B}^d}_\alpha \right) +\norm{\nabla \sigma}_\infty^2 \norm{{\bf B}^d}_\alpha^2  |\sigma(0)| \right)   \right. \\
		& \left. +   |\sigma(0)|  \left(\norm{{\bf B}^d}_\alpha + \norm{{\bf B}^d}_\alpha\norm{\nabla \sigma}_\infty  \mu  +  {\mu^{\frac{1}{\alpha}-1}}  \right) \right) \\
		C_{16} (\mu) &:= C_{13}\norm{\nabla \sigma}_\infty  \left( C_{14}\mu^{2} \norm{\nabla \sigma}_\infty^2 \norm{{\bf B}^d}_\alpha^2    + \norm{{\bf B}^d}_\alpha + \mu \norm{\nabla \sigma}_\infty    \norm{{\bf B}^d}_\alpha  + {\mu^{\frac{1}{\alpha}-1}}    \right) \\
		C_{13} &:=  \frac{1}{(1-{3} M)(1 - K_2)} \\
		C_{14} &:= \frac{K}{1 -K_2}
		\end{align*}	
		for a small enough time horizon $\mu$ satisfying 
		\begin{subequations}
			\label{eq:muconds}
			\begin{align}
			\label{eq:mucond1}
			\norm{{\bf B}^d}_\alpha  \norm{\nabla \Sigma}_\infty  \mu^{2} &\leq M \\
			\label{eq:mucond2}
			 \norm{{\bf B}^d}_\alpha \norm{\nabla \sigma}_\infty \mu&\leq M \\
			\label{eq:mucond3}
			{ \norm{\nabla b}_\infty \mu^{\frac{1}{\alpha}}} &\leq M \\
			\label{eq:mucond4}
			 \norm{\nabla \sigma}_\infty \norm{{\bf B}^d}_\alpha \mu &\leq \frac{K_2}{K}  \\
			\label{eq:mucond5}
			\mu^{2} (C_4(\mu) + C_3(\mu)) &\leq \frac{K_2 (1-3M)(1-K_2)}{K}
			\end{align}
			with
			\begin{align*}
			C_3(\mu) &:= \norm{\nabla \Sigma}_\infty \left((1 + \mu^{2})\norm{{\bf B}^d}_\alpha^2 + \norm{{\bf B}^d}_\alpha \right) \\
			C_4(\mu) &:= \norm{\nabla \sigma}_\infty^2 \norm{{\bf B}^d}_\alpha^2 \mu 
			\end{align*}
		\end{subequations}
		and where $0 < M < \frac{1}{3}$, $0 < K_2 < 1$ and $K = f(\theta)$.  
	\end{lem}
	
\begin{lem}
	\label{lem:uniformest_sde}
	\textbf{Uniform estimates on the smooth approximation}
	Under Assumptions \ref{ass:bsig} and \ref{ass:correc}, and when $x_0$ is a fixed value, it holds that 
	\begin{align*}
	\sup_{d \in \mathbb{N}} \mathbb{E}\left[ \sup_{t \in [0,T]} |X_t^d|^2  \right] < \infty
	\end{align*}

	\begin{proof}
		Before extending to the large time horizon, we require a precise expression for $K^*$ based on the conditions on $\mu$ in (\ref{eq:muconds}).  In particular, 
		\begin{align}
		\label{eq:condmu1_4}
		\mu \leq \min \left( \frac{1}{\norm{{\bf B}^d}_\alpha^{1/2}}, \frac{1}{\norm{{\bf B}^d}_\alpha}, 1  \right) C(M, \norm{\nabla \Sigma}_\infty , \norm{\nabla b}_\infty, \norm{\nabla \sigma}_\infty, \alpha)
		\end{align}   
		ensures $\mu$ satisfies (\ref{eq:mucond1})-(\ref{eq:mucond4}).  To analyse the final condition (\ref{eq:mucond5}), let $K_3:=  \frac{K_2 (1 - K_2)}{K} (1 - 3M)$ and $a:= \norm{\nabla \sigma}_\infty$ and $\beta:=\norm{\nabla \Sigma}_\infty$.  We then require 
		\begin{align}
		\label{eq:longcond}
		\mu^{3}a^2 \norm{{\bf B}^d}_\alpha^2 + \beta \left(\norm{{\bf B}^d}_\alpha^2 + \norm{{\bf B}^d}_\alpha \right) \mu^{2} + \beta \norm{{\bf B}^d}_\alpha^2 \mu^{4} \leq K_3.
		\end{align}
		It is not difficult to see that (\ref{eq:longcond}) holds in the case $\mu \geq 1$ when 
		\begin{align*}
		\left(a^2 \norm{{\bf B}^d}_\alpha^2 + \beta  \left(2 \norm{{\bf B}^d}_\alpha^2 + \norm{{\bf B}^d}_\alpha \right) \right)\mu^{4} \leq K_3
		\end{align*}
		and in the case $\mu < 1$ when 
		\begin{align*}
		\left(a^2 \norm{{\bf B}^d}_\alpha^2 + \beta  \left(2 \norm{{\bf B}^d}_\alpha^2 + \norm{{\bf B}^d}_\alpha \right) \right)\mu^{2} \leq K_3.
		\end{align*}
		Together, these conditions imply that we require 
		\begin{align*}
		\mu &\leq \min \left( \left(\frac{K_3}{(a^2 + 2\beta ) \norm{{\bf B}^d}_\alpha^2 + \norm{{\bf B}^d}}_\alpha \right)^{1/2},  \left(\frac{K_3}{(a^2 + 2\beta ) \norm{{\bf B}^d}_\alpha^2 + \norm{{\bf B}^d}}_\alpha \right)^{1/4}  \right). 
		\end{align*}
		Combining the above with (\ref{eq:condmu1_4}) means that (\ref{eq:mucond1})-(\ref{eq:mucond5}) are satisfied when 
		\begin{align*}
		\mu & \leq \min \left( \frac{1}{ \left[(a^2 + 2\beta ) \norm{{\bf B}^d}_\alpha^2 + \norm{{\bf B}^d}_\alpha \right]^{1/2}},  \frac{1}{ \left[(a^2 + 2\beta ) \norm{{\bf B}^d}_\alpha^2 + \norm{{\bf B}^d}_\alpha \right]^{1/4}},  \frac{1}{\norm{{\bf B}^d}_\alpha^{1/2}}, \frac{1}{\norm{{\bf B}^d}_\alpha}, 1   \right) \\
		& \times C(M, \norm{\nabla \Sigma}_\infty, \norm{\nabla b}_\infty, \norm{\nabla \sigma}_\infty, K, K_2, \alpha).
		\end{align*}
		Again, considering the cases $\norm{{\bf B}^d}_\alpha \leq 1$ and $\norm{{\bf B}^d}_\alpha > 1$ separately we obtain the following expression 
		\begin{align}
		\label{eq:Kstarbound}
		K^* &= \frac{1}{1 + \norm{{\bf B}^d}}_\alpha C^*(M, \norm{\nabla \Sigma}_\infty, \norm{\nabla b}_\infty, \norm{\nabla \sigma}_\infty, K, K_2, \alpha).
		\end{align}	
		Furthermore, (\ref{eq:mucond2}) and (\ref{eq:mucond3}) imply that 
		\begin{align*}
		C_{16} (\mu) &\leq C_{13} \norm{\nabla \sigma}_\infty  \left( \norm{{\bf B}^d}_\alpha + C_{14}M^2 + M +  {C_{17}}    \right) \\
		&:= C_{18} \norm{{\bf B}^d}_\alpha + C_{19}
		\end{align*}
		where $C_{17}:= \left(  \frac{M}{\norm{\nabla b}_\infty}  \right)^{1 - \alpha}$.  
		Likewise, (\ref{eq:mucond1}), (\ref{eq:mucond2}) and (\ref{eq:mucond3}) imply 
		\begin{align*}
		C_{15}(\mu) & \leq C_{13} \left( ( C_{14}M  + |\sigma(0)|)  \norm{{\bf B}^d}_\alpha + 1 + C_{14}M \left(  1 + \frac{M}{\norm{\nabla \Sigma}_\infty}  \right)   +  C_{14} |\sigma(0)| M^2  +  |\sigma(0)|  ( M  +  {C_{17}}  )  \right) \\
		&:= C_{20} \norm{{\bf B}^d}_\alpha + C_{21}.
		\end{align*}
		Therefore, when $\mu \leq K^*$, we have that 
		\begin{align*}
		\norm{\hat{X}^d}_\alpha  \leq C_{20} \norm{{\bf B}^d}_\alpha + C_{21}  +  |x_0| \left( C_{18} \norm{{\bf B}^d}_\alpha + C_{19}   \right). 
		\end{align*}
		Then by Proposition 7 in \cite{Lejayglob2012},
		\begin{align}
		\label{eq:supybound}
		\sup_{t \in [0,T]} |\hat{X}^d_t| & \leq R(T)|x_0| + R(T)\frac{C_{20} \norm{{\bf B}^d}_\alpha + C_{21}  }{C_{18} \norm{{\bf B}^d}_\alpha + C_{19}}, 
		\end{align}
		with 
		\begin{align*}
		R(T) &= \exp \left(C_{16}(\mu)\left(1 + \frac{1}{K^*} \right)^{1 - \alpha} \right) \exp(\max \{ T, \mu  \}) \\
		& \leq \exp \left(C_{16}(\mu)\left(1 + \frac{1}{K^*} \right) \right) \exp(\max \{ T, \mu  \}),
		\end{align*}
		since $\frac{1}{3} < \alpha < \frac{1}{2}$.  Furthermore, 
		\begin{align*}
		\frac{C_{20} \norm{{\bf B}^d}_\alpha + C_{21}}{C_{18}\norm{{\bf B}^d}_\alpha + C_{19}} & \leq \frac{C_{20}}{C_{18}} + \frac{C_{21}}{C_{19}} \\
		& := C_{22}
		\end{align*}
		and also using (\ref{eq:Kstarbound}) we have 
		\begin{align*}
		\exp \left(C_{16}(\mu) \left(1 + \frac{1}{K^*} \right) \right) & \leq \exp(C_2^*C_{19})\cdot \exp\left(C_2^*   C_{18} \norm{{\bf B}^d}_\alpha \right), 
		\end{align*}
		where $C_2^* := \frac{C^* +1}{C^*}$.  Finally, we can conclude using (\ref{eq:supybound}) that 
		\begin{align*}
		\sup_{t \in [0,T]} |\hat{X}^d_t|^2 & \leq \exp \left(2C_2^*C_{19} \right) \exp(2\max \{ T, \mu  \})  \left(|x_0| + C_{22} \right)^2 \exp \left(2C_2^*C_{18} \norm{{\bf B}^d}_\alpha \right)  \\
		& \leq \exp \left(2C_2^*C_{19} + 2\max \{ T, \mu  \} + 2C_2^*C_{18} \norm{{\bf B}^d}_\alpha \right)  \left(|x_0| + C_{22} \right)^2 \\
		& := C_{24} \exp(2C_2^*C_{18} \norm{{\bf B}^d}_\alpha)
		\end{align*}
		where $C_{2}^*, C_{18}, C_{19}, C_{22}$ are positive constants depending only on $x_0, \sigma(0), M, \norm{\nabla \Sigma}_\infty, \norm{\nabla b}_\infty, \norm{\nabla \sigma}_\infty, K, K_2, \alpha$.  Finally, for nested piecewise linear approximations, Theorem 13.19 in \cite{Friz2009} gives us that 
		\begin{align*}
		\exp(C_2^*C_{18} \norm{{\bf B}^d}_\alpha) \leq \exp(C_2^*C_{18} M_g) < \infty \quad \forall \enskip d \in \mathbb{N}
		\end{align*} 
		where $M_g$ is a positive random variable with Gaussian tails independent of $d$ and $M_g < \infty$ almost surely. Taking expectation of both sides gives the desired result.
		
	\end{proof}
	
\end{lem}

	We are now ready to show that solutions of localised ODEs converge to the solution of the ODE of interest, uniformly in the approximation parameter. 

	\begin{lem}
		\label{lem:uniform}
		\textbf{Uniform convergence in $d$.}  For every $T > 0$ and $X_t^{d,n}$ satisfying (\ref{eq:sden}) and $X_t^d$ satisfying (\ref{eq:smoothsde}), 
		\begin{align}
		\label{eq:unifconv}
		\lim_{n \rightarrow \infty} \sup_{d \in \mathbb{N}} \mathbb{E} \left[\sup_{0 \leq t \leq T} |X_t^{d,n} - X_t^d|^2 \right] = 0 
		\end{align}
		\begin{proof}
			For each fixed $d$, define 
			\begin{align*}
			\tau^d_n := \inf \{ t \geq \delta_d \enskip | \enskip |X_{t-\delta_d}^{d,n} |^2 \notin S_n \cup |X_{t-\delta_d}^d |^2 \notin S_n  \} - \delta_d 
			\end{align*}  
			 which is a stopping time with respect to $\mathcal{F}_t$.  Define
			\begin{align*}
			\mathcal{S}^{dn}_T := \sup_{0 \leq t \leq T} |X_t^{d,n} - X_t^d|^2. 
			\end{align*}  
			It holds that
			\begin{align*}
			\mathbb{E}[\mathcal{S}^{dn}_T] = \mathbb{E}[\mathds{1}_{\tau_n^d > T} \mathcal{S}^{dn}_T] + \mathbb{E}[\mathds{1}_{\tau_n^d \leq T} \mathcal{S}^{dn}_T] 
			\end{align*}       
			and by Hoelder inequality, 
			\begin{align*}
			\mathbb{E}[\mathds{1}_{\tau_n^d \leq T} \mathcal{S}^{dn}_T] & \leq \mathbb{E}[(\mathds{1}_{\tau_n^d \leq T})^p]^{1/p} \mathbb{E}[(\mathcal{S}^{dn}_T)^q]^{1/q} \\
			& = (\mathbb{P}(\tau_n^d \leq T))^{1/p}  \mathbb{E}[(\mathcal{S}^{dn}_T)^q]^{1/q}, 
			\end{align*} 
			where $\frac{1}{q} + \frac{1}{p} = 1$.  For the case $\tau_n^d > T$, we start by considering the Stratonovich step 2 lift of $B^d_t$ as in Lemma \ref{lem:uniformest_sde}.  We can again work with the RDEs whose solutions coincide with the ODEs.  Continuity of the Ito-Lyons map under Assumptions \ref{ass:bsig} and \ref{ass:correc} is established in Theorem 1 in \cite{Lejayglob2012}, which implies that for all $t < \tau_n^d$,
			\begin{align*}
			|X_t^{d,n} - X_t^d| = 0 \quad \text{a.s.}
			\end{align*}
			since both $X_t^{d,n}$ and $X_t^d$ are driven by the same ${\bf B}_t^d$ and $X_0^{d,n} = X_0^d$.  
			Furthermore, it follows directly that $\mathbb{E}[\mathcal{S}^{dn}_T] = 0$ when $\tau_n^d > T$ for every fixed $n,d$ (and uniformly in $d$).  Therefore, we have that for fixed $n$, 
			\begin{align*}
			\sup_{d \in \mathbb{N}} \mathbb{E}[\mathcal{S}^{dn}_T] & \leq \sup_{d \in \mathbb{N}} (\mathbb{P}(\tau_n^d \leq T) )^{1/p} \mathbb{E}[(\mathcal{S}^{dn}_T)^q]^{1/q}. 
			\end{align*}
			Lemma \ref{lem:uniformest_sde} clearly holds for $X_t^{d,n}$ uniformly in $n$, therefore we have that 
			\begin{align*}
			\mathbb{E}[(\mathcal{S}_{dn})^q]^{1/q} & \leq \left(\mathbb{E} \left[ \sup_{0 \leq t \leq T} |X_t^{d,n}|^{2q} \right] + \mathbb{E} \left[ \sup_{0 \leq t \leq T} |X_t^{d}|^{2q} \right] \right)^{1/q} \\
			& < \infty 
			\end{align*}
			uniformly in $n$ and $d$.  Furthermore, by standard arguments we have that 
			\begin{align*}
			\mathbb{P}(\tau_n^d \leq T) &= \mathbb{P} \left(  \sup_{t \in [0, T]} |X^d_t|^2 > n \right) \\
			& = \mathbb{E} [\mathds{1}_{\{  \sup_{t \in [0, T]} |X^d_t|^2 > n  \}}] \\
			& \leq \frac{1}{n} \mathbb{E} \left[\sup_{t \in [0, T]} |X^d_t|^2 \right]
			\end{align*}
			Then taking sup wrt $d$ on both sides in the above, and using $\sup_{d \in \mathbb{N}} \mathbb{E} \left[\sup_{t \in [0, T]} |X^d_t|^2 \right] < \infty$ from Lemma \ref{lem:uniformest_sde} implies
			\begin{align*}
			\lim_{n \rightarrow \infty} \sup_{d \in \mathbb{N}} (\mathbb{P}(\tau_n^d \leq T))^{1/p}  = 0. 
			\end{align*} 
			Combining all gives 
			\begin{align*}
			\lim_{n \rightarrow \infty} \sup_{d \in \mathbb{N}} \mathbb{E}[\mathcal{S}_{dn}]  &\leq  \left(\lim_{n \rightarrow \infty} \sup_{d \in \mathbb{N}} (\mathbb{P}(\tau_n^d \leq T))^{1/p} \right)  \left( \lim_{n \rightarrow \infty} \sup_{d \in \mathbb{N}} \mathbb{E}[(\mathcal{S}_{dn})^q]^{1/q}  \right)  \\
			& = 0.
			\end{align*}
			since both limits are finite. 
			 
		\end{proof}
		
	\end{lem}

	The final ingredient is the pointwise convergence statement, as stated in the following lemma. 
	\begin{lem}
		\label{lem:pointwise}
		\textbf{Pointwise convergence in $n$.}  For every $T > 0$ and $X_t^{d,n}$ satisfying (\ref{eq:smoothlocapp}) and $X_t^n$ satisfying (\ref{eq:sden}), 
		\begin{align}
		\label{eq:pointconv}
		\lim_{d \rightarrow \infty} \mathbb{E} \left[\sup_{0 \leq t \leq T} |X_t^{d,n} - X_t^n|^2 \right] = 0  \quad \forall \enskip n \in \mathbb{N} 
		\end{align}
		\begin{proof}
			Follows trivially from Theorem 7.2 in \cite{Ikeda1989} since for every fixed $n$, $b^n$ and $\sigma^n$ are uniformly bounded. 
		\end{proof}
	\end{lem}
	 
	Since $\delta_2(X,Y) := \mathbb{E} \left[\sup_{0 \leq t \leq T} |X_t - Y_t|^2 \right]$ is a distance metric on the space of $L^2$ integrable processes\footnote{Here we are working with equivalence classes}, to which all proceses considered here belong, we can apply the Moore-Osgood Theorem (see Theorem \ref{theo:mooreosgood}) together with Lemmas \ref{lem:pointwise}, \ref{lem:uniform} and \ref{lem:locsde} to obtain 
	\begin{align*}
		\lim_{n \rightarrow \infty} \lim_{d \rightarrow \infty} \mathbb{E} \left[\sup_{0 \leq t \leq T} |X_t^{d,n} - X_t^n|^2  + \sup_{0 \leq t \leq T} |X_t^{n} - X_t|^2 \right]& = \lim_{d \rightarrow \infty} \lim_{n \rightarrow \infty} \mathbb{E} \left[\sup_{0 \leq t \leq T} |X_t^{d,n} - X_t^d|^2  + \sup_{0 \leq t \leq T} |X_t^{d} - X_t|^2 \right] \\
		& = 0.
	\end{align*}
	Furthermore, since uniform convergence implies pointwise convergence, we have from Lemma \ref{lem:uniform} that for any fixed $d \in \mathbb{N}$, 
	\begin{align*}
	\lim_{n \rightarrow \infty} \mathbb{E} \left[\sup_{0 \leq t \leq T} |X_t^{d,n} - X_t^d|^2  \right] = 0 
	\end{align*}
	which then gives the desired result.  This concludes the proof of Theorem \ref{theo: sdeconv}. 
	
	\newpage 
		 
	 \section*{Acknowledgements}
	 This research has been partially funded by Deutsche Forschungsgemeinschaft (DFG)- SFB1294/1 - 318763901.  The author is grateful to Wilhelm Stannat and Sebastian Reich for helpful feedback on this work.

 	 \appendix
	 \section{Appendix}	
	
	 \begin{theorem}
	 	\label{theo:mooreosgood}
	 	\textbf{Moore-Osgood Theorem.} Let $(M, d_M)$ be a metric space and $(m_{d,n})_{d,n \in \overline{\mathbb{N}}}$ be a sequence in $M$ where $\overline{\mathbb{N}} = \mathbb{N} \cup \infty$.  If 
	 	\begin{itemize}
	 		\item [(i)] $\lim_{n \rightarrow \infty} \sup_{d \in \mathbb{N}} d_M (m_{d,n}, m_{d, \infty}) = 0$ and 
	 		\item [(ii)] $\lim_{d \rightarrow \infty} d_M(m_{d,n}, m_{\infty, n}) = 0$ for all $n \in \mathbb{N}$,
	 	\end{itemize}
	 	then the joint limit $\lim_{d,n \rightarrow \infty} m_{n,d}$ exists.  In particular, it holds that 
	 	\begin{align*}
	 	\lim_{d,n \rightarrow \infty} m_{d,n} = \lim_{n \rightarrow \infty} m_{\infty, n} = \lim_{d \rightarrow \infty} m_{d, \infty}.
	 	\end{align*}
	 \end{theorem}

	 \subsection*{Proof of Lemma \ref{lem:ybound}}
	 \label{sec:proofrdebound}

	 		Define 
	 		\begin{align*}
	 		D(s,t):= \int_s^t b(\hat{X}^d_u) du - \sigma(\hat{X}^d_s){\bf B}_{s,t}^{d,1} - \Sigma(\hat{X}^d_s){\bf B}_{s,t}^{d,2}. 
	 		\end{align*}	
	 		
	 		Using the fact that $b$ is Lipschitz, we have that 
	 		\begin{align*}
	 		\left| \int_s^t b(\hat{X}_u^d) du \right| & \leq \int_s^t |b(\hat{X}_u^d) - b(x_0)| du + \int_s^t |b(x_0)| du  \\
	 		& \leq \norm{\nabla b}_\infty \int_s^t |\hat{X}_u^d - x_0| du + |b(x_0)| (t-s) \\
	 		& \leq \norm{\nabla b}_\infty \norm{\hat{X}^d}_\alpha  \int_s^t  u^{\alpha} du + |b(x_0)| (t-s) \\
	 		& \leq (\norm{\nabla b}_\infty \norm{\hat{X}^d}_\alpha  \mu  + |b(x_0)| )(t-s)^{\alpha} \mu^{\frac{1}{\alpha} -1}
	 		\end{align*}
	 		Similarly, using the fact that $\sigma$ is Lipschitz, 
	 		\begin{align*}
	 		|\sigma(\hat{X}_s^d)| |{\bf B}_{s,t}^{d,1}| & \leq \norm{\nabla \sigma}_\infty \frac{|\hat{X}^d_s - x_0|}{s^{\alpha}} s^{\alpha} \norm{{\bf B}^d}_\alpha (t-s)^{\alpha} + |\sigma(x_0)| \norm{{\bf B}^d}_\alpha (t-s)^{\alpha}  \\
	 		& \leq \norm{\nabla \sigma}_\infty \norm{\hat{X}^d}_\alpha  \mu \norm{{\bf B}^d}_\alpha(t-s)^{\alpha} + |\sigma(x_0)|\norm{{\bf B}^d}_\alpha (t-s)^{\alpha}
	 		\end{align*}
	 		 By similar arguments we have 
	 		\begin{align*}
	 		|\Sigma(\hat{X}_s^d)| |{\bf B}_{s,t}^{d,2}|  & \leq  |\Sigma(\hat{X}^d_s) - \Sigma(x_0)| \frac{|{\bf B}_{s,t}^{d,2}|}{(t-s)^{2\alpha}} (t-s)^{2\alpha} + |\Sigma(x_0)| \frac{|{\bf B}_{s,t}^{d,2}|}{(t-s)^{2\alpha}} (t-s)^{2\alpha} \\
	 		& \leq \norm{\nabla \Sigma}_\infty  \norm{\hat{X}^d}_\alpha \mu^{2} \norm{{\bf B}^d}_\alpha (t-s)^{\alpha} + |\sigma(x_0)| \norm{\nabla \sigma}_\infty \norm{{\bf B}^d}_\alpha (t-s)^{\alpha} \mu 
	 		\end{align*}
	 		Then using the existence of a solution in the sense of Davie (\ref{eq:solndavie}) and the above estimates, we have that 
	 		\begin{align*}
	 		\norm{\hat{X}^d}_\alpha \leq L \mu^{2}+ { C_{61}(\mu)\norm{\hat{X}^d}_\alpha + C_{62}(\mu, x_0)} + C_7(\mu, x_0) + C_8(\mu) \norm{\hat{X}^d}_\alpha + C_9(\mu) \norm{\hat{X}^d}_\alpha 
	 		\end{align*}
	 		where 
	 		\begin{align*}
	 		C_7(\mu, x_0) &:=  \norm{{\bf B}^d}_\alpha|\sigma(x_0)|  (1+ \norm{\nabla \sigma}_\infty  \mu) \\
	 		C_8(\mu) & := \norm{{\bf B}^d}_\alpha  \norm{\nabla \Sigma}_\infty  \mu^{2} \\
	 		C_9(\mu) &:= \norm{{\bf B}^d}_\alpha \norm{\nabla \sigma}_\infty \mu \\
	 		C_{61}(\mu) &:= { \norm{\nabla b}_\infty \mu^{\frac{1}{\alpha}}} \\
	 		C_{62}(\mu, x_0) &:= {|b(x_0)| \mu^{\frac{1}{\alpha}-1}}.
	 		\end{align*}
	 		If we consider a short time horizon such that 
	 		\begin{align}
	 		\label{eq:conds1}
	 		C_8(\mu) \leq M, \quad C_9(\mu) \leq M, \quad {C_{61}(\mu) \leq M} 
	 		\end{align}
	 		we have that for the case $M < \frac{1}{3}$, 
	 		\begin{align}
	 		\label{eq:intermybound}
	 		\norm{X^d}_\alpha \leq \frac{1}{1-3M} \left(   L\mu^{2} +  C_{62}(\mu, x_0) + C_7(\mu, x_0)   \right).
	 		\end{align}
	 		Consider 
	 		\begin{align*}
	 		D(s,r,t):= D(s,t) - D(s,r) - D(r,t)
	 		\end{align*}
	 		which by Lemma 1 in \cite{Lejayglob2012} is an almost additive functional and can be bounded using similar arguments as developed above for $D(s,t)$.
	 		The sewing lemma (see e.g. Lemma 4.1 in \cite{Coutin2012}) implies that there exists a universal constant $K$ depending only on $\theta$ such that
	 		\begin{align}
	 		\label{eq:Lprime}
	 		L \leq  \underbrace{\sup_{0\leq s < t \leq T} \frac{|\hat{X}^d_t - \hat{X}^d_s - D(s,t)|}{(t-s)^{\theta}}}_{:=L'} \leq K \sup_{0 \leq s < r< t \leq T} \frac{|D(s,r,t)|}{(t-s)^\theta}
	 		\end{align}
	 		and by the same arguments as in the proof of Proposition 2 in \cite{Lejayglob2012}, which are unaffected by the presence of the drift in the RDE, we have for the small time horizon satisfying 
	 		\begin{align}
	 		\label{eq:cond2}
	 		KC_5(\mu) \leq K_2 
	 		\end{align}
	 		that
	 		\begin{align*}
	 		L' \leq \frac{K}{1 -K_2} (C_{10}(\mu) + (C_4(\mu) + C_3(\mu)) \norm{\hat{X}^d}_\alpha + C_{12}(\mu) |x_0|) 
	 		\end{align*}
	 		where $0 < K_2 < 1$ and 
	 		\begin{align*}
	 		C_3(\mu) &:= \norm{\nabla \Sigma}_\infty ((1 + \mu^{2})\norm{{\bf B}^d}_\alpha^2 + \norm{{\bf B}^d}_\alpha) \\
	 		C_4(\mu) &:= \norm{\nabla \sigma}_\infty^2 \norm{{\bf B}^d}_\alpha^2 \mu \\
	 		C_5(\mu) &:= \norm{\nabla \sigma}_\infty \norm{{\bf B}^d}_\alpha \mu \\
	 		C_{11}(\mu) &:=  \norm{\nabla \sigma}_\infty^2 \norm{{\bf B}^d}_\alpha^2  \\
	 		C_{12}(\mu) &:= C_{11} \norm{\nabla \sigma}_\infty \\
	 		C_{10} (\mu) & := C_3(\mu) + |\sigma(0)| C_{11}(\mu).
	 		\end{align*}
	 		Substituting into (\ref{eq:intermybound}) gives
	 		\begin{align*}
	 		\norm{X^d}_\alpha \leq \frac{1}{1-{3} M} \left(  \frac{K}{1 -K_2}\mu^{2} \left(C_{10}(\mu) + (C_4(\mu) + C_3(\mu))  \norm{X^d}_\alpha + C_{12}(\mu) |x_0| \right) + C_7(\mu, x_0) + { C_{62}(\mu, x_0)}   \right)
	 		\end{align*}
	 		Again restricting to a short time horizon, we have that if for example 
	 		\begin{align}
	 		\label{eq:cond}
	 		\frac{1}{1-{3} M}  \frac{K}{1 -K_2}\mu^{2} (C_4(\mu) + C_3(\mu)) \leq K_2
	 		\end{align}
	 		then (\ref{eq:yroughbound}) holds. Condition (\ref{eq:muconds}) easily follows from combining (\ref{eq:conds1}), (\ref{eq:cond2}), (\ref{eq:cond}).

\bibliographystyle{alpha}
\bibliography{literature_WZ.bib}

\end{document}